\documentclass [12pt]
{amsart}

\author{Saharon Shelah}
\title[ccc with Sacks/Laver]{Consistently there is no non trivial ccc   
 forcing notion with the Sacks or Laver  property   
}

\usepackage{amsthm}
\usepackage{amssymb} 
\usepackage{a4wide}
\swapnumbers

\newcounter{conditions}  
\def\itm#1\par{\addtocounter{conditions}{1}
\begin{itemize}
\item[$(*)_{\arabic{conditions}}$] #1
\end{itemize}
}

\makeatletter
 \def\LABEL#1{\def\@currentlabel{$(*)_{\arabic{conditions}}$}\label{#1}}

\def\on{{\upharpoonright}}

\newcommand{\M}{{\mathbb M}}
\newcommand{\R}{{\mathbb R}}
\newcommand{\tlo}{ {}^{\omega>} 2}
\newcommand{\forces}{\Vdash}
\newcommand{\Q}{{\mathbb Q}}
\renewcommand{\P}{{\mathbb P}}
\newcommand{\T}{{\mathbb T}}
\newcommand{\dom}{{\rm Dom}}
\newcommand{\Dom}{{\rm Dom}}
        
\def\itema$#1${\item[$#1$]}

 \def\name#1{\mathpalette\putsimunder{#1}}
\def\putsimunder#1#2{\oalign{$#1#2$\crcr\hidewidth
\vbox to.2ex{\hbox{$#1\widetilde{\hphantom{#2}}$}\vss}\hidewidth}}

\dedicatory{Dedicated to the memory of Paul Erd\H os} 
\thanks{The author thanks the Israel Science Foundation (founded by
the Israel  
 Academy of Sciences) for partially 
 supporting this research. Pub number 723}

\def\xx#1 {\newtheorem{#1}[thm]{#1}}
\xx Theorem
\xx Question
\xx Lemma

\theoremstyle{definition}
\xx Fact
\xx Framework
\xx Remark
\xx Notation
\xx Claim
\xx Observation
\xx Definition

\begin{document}

\maketitle

\section{Introduction}

At the recent  set theory conference, Boban 
 Velickovic asked the following question: 

\begin{Question} Is there a nontrivial forcing notion with the
 Sacks property which is also ccc? 
\end{Question}
(See below for a definition of the Sacks property.)

 A ``definable'' variant  of this question 
 has been answered in \cite{Sh:480}: 
\begin{quote} Every  nontrivial Souslin forcing notion which has the
 Sacks property has an uncountable antichain. 
\end{quote}
(A Souslin forcing notion is a forcing notion for which the set of
 conditions, the comparability relation and the incompatibility
 relation are all analytic subsets of the reals.  See \cite{JdSh:292}
 and   \cite{Sh:480} for details). 

We show here

\begin{Theorem}\label{jlx}
The following statement is equiconsistent with ZFC: 
\begin{itemize}
\item[$(*)$] Every  nontrivial forcing notion which has the
 Sacks property has an uncountable antichain. 
\end{itemize}
\end{Theorem}

Our proof follows the ideas from \cite{Sh:480}. 

Independently, Velickovic  
 has also proved the  consistency of $(*)$,  following 
\cite{Sh:480}   and some of his works.   In fact, he shows
% \begin{Theorem}\label{vel} The 
that
 the proper forcing axiom (PFA),   and even the  open
 coloring axiom implies $(*)$. 
% \end{Theorem}

Our proof shows that also the following strengthening of $(*)$: 
\begin{itemize}
\item[$(**)$] Every  nontrivial forcing notion which has the
 Laver  property has an uncountable antichain. 
\end{itemize}
is equiconsistent with ZFC.

 Note that if cov(meagre) $=$ continuum  
 (which follows e.g.\ from PFA) then there  
 is a (non principal) Ramsey ultrafilter  
 on  
 $\omega$.   The ``Mathias'' forcing notion for shooting making this
 ultrafilter principal has the Laver
 property and is ccc, so $(**)$ does not follow from PFA. 

 So our result  and Velickovic' result are incomparable.

\begin{Definition}
Let $g\in \omega ^ \omega$ be increasing. A
{\em  $g $-slalom}
is a sequence $\bar A=(A_n:n\in \omega )$, $A_n \subseteq \omega $, 
$|A_n|\le g(n)$.  We say that $\bar A$ covers $\eta\in \omega^\omega$ 
iff $ \forall n\, \eta(n)\in A_n$.

Let $V_1 \subseteq V_2 $ e models of set theory.  We say that
$(V_1,V_2)$ is {\em $(f,g)$-bounding} iff: 
\begin{quote}
For all $\eta\in \prod_n f(n) \cap V_2$ there is a $g$-slalom in $V_1$
covering $f$.
\end{quote}
We say that a forcing notion $\R$ is $(f,g)$-bounding if the pair $(V,
V^\R)$ is  $(f,g)$-bounding 

A pair $(V_1, V_2)$ has the {\em Laver property} iff $(V_1, V_2)$ is
$(f,g)$-bounding for all increasing $f,g\in V_1$, or in other words:
For all $g\in V_1$, every function in $V_2$ which is bounded by a
function in $V_1$ is covered by a $g$-slalom from $V_1$. 

Similarly, $\R$  has the Laver property iff $(V, V^\R)$ has the Laver
property.

$(V_1,V_2)$ has the Sacks property if it has the Laver property, and
every function in $V_2$ is bounded by a function in $V_1$.

\end{Definition}

\section{A lemma on Mathias forcing}

The following lemma is a theorem of ZFC. 

\begin{Lemma}\label{mainlemma}
Let $\M$ be the Mathias forcing, $\name\eta$ the increasing 
enumeration of the
$\M$-generic subset of $\omega $, and assume that $\name T$ is  a name 
of a tree $ \subseteq \tlo $ such that  
$$ (!) \qquad \qquad 
 \forces _\M  | \name T \cap {}^ {\name\eta(n)}2 | \le
2^{\name\eta(n-1)}$$
Then there is a family $(T_i, q_i: i< 2^{\aleph_0} ) $ such that  
for all $i,j$:
\begin{enumerate}
\item $ q_i\in \M$
\item $T_i \subseteq \tlo $ is a tree
\item $q_i \forces \name T \subseteq T_i$
\item Whenever $i\not=j$, then $\lim T_i\cap \lim T_j$ is finite. 
\\
(We write $\lim T$ for the set of branches of a tree $T$.)
\end{enumerate}
\end{Lemma}

\begin{proof}

Conditions in the Mathias forcing are of the form $(w,A)$, where 
 $w$ is a  
 finite set of natural numbers 
  and $A$  is an infinite subset 
  of $\omega$ 
 such that   $\sup(w ) < \min A $.   If $(w,A)\in \M$, and 
$w = \{ w_0< \cdots < w_{\ell-1}\}$, then 
$$ (w,A) \forces \name\eta(i)=w_i \ (i= 0,\ldots,\ell-1) , \qquad 
      \name \eta(i)\in A  \ (i=\ell, \ell+1, \ldots) $$

Using the fact that truth values in $V^\M$ can be
decided by pure extensions, we 
can find a condition $(w^*, A^*)$ which forces that $\name T$ is
``decided continuously'' by $\name \eta$, more specifically: 
\itm 	\LABEL{continuous}
there is a function $t$ with domain $[A^*]^{<\omega }$,
 $ \forall w:  t(w) \subseteq {}^{\omega > } 2$,  such that  
\begin{quote}
For all finite  $u \subseteq A^*$ and all $n$ with
 sup$(u) <n \in  A^*$: \\
 $(w^*\cup u,  A^* \setminus (n+1)) \forces 
          \name T_{\delta}\cap {}^{n} 2 = t(w) \cap {}^n 2$. 
\end{quote}

 For 
  any finite $u \subseteq A^*$ let  
 $$T_{u} =^{df} \{ \eta \in {}^{\omega>} 2 :
\mbox{ for 
  any large enough $k<\omega$ we have 
 $(w^* , A^* \setminus k ) \Vdash_{{\Bbb Q}}$ ``$\eta 
 \in {\name T}$''}  \}$$
So $T_{u}$ 
 is a subtree of {}$^{\omega>}2$.

We have
\itm If $w^*\cup u = \{ w_0< \cdots < w_{\ell-1}\}$, then 
$ (w^*\cup u, A \setminus k)  \forces \name \eta(\ell-1) = \max u$,
\ $ \name \eta(\ell ) \ge k $\\
hence (recall condition (!)): 
$ (w^*\cup u, A \setminus k)  \forces | \name T \cap {}^k 2 |
\le  2^ {\eta(\ell-1)} = 2^{ \max u} $. 

It follows that: 
\itm
 for every finite $u \subseteq A^*$,  
 non empty for simplicity, we have 
 $$k<\omega \rightarrow \vert T_{u} \cap  {}^{k} 2 \vert \le  
2^{{\rm max} (u ) } , $$ hence lim$(T_{u})$  
 is a finite subset of ${}^\omega 2$. 

We also get: 
\itm   if $u \cup \{m,k\} \subseteq A^*$  
 and sup$(u) <m < k$,
 then $T_{u} \cap {}^{m} 2 =$ 
  $T_{u \cup \{k\}}  \cap {}^{m} 2$. 

%  and $T_{u} \cap {}^{{\rm sup}(u )} 2 =$  
%  $T_{u \cup \{m\}}$  
%   $\cap {}^{{\rm sup}(u )} 2$ 

[Proof: 
We know that 
already $ (w^*\cup u, A \setminus m)$ decides $\name T\cap {}^m 2  $,
 and the conditions 
$ (w^*\cup u\cup \{m\}, A \setminus k)$ and 
 $ (w^*\cup u, A \setminus k)$  (for $k>m$) are both
stronger than  $ (w^*\cup u, A \setminus m)$.]

In particular, we get: 
\itm  for all  finite  $u \subseteq A^*$ the 
  sequence  
 $\langle  T_{u \cup  \{ k \} } : k \in$ 
  $A^* \rangle$ converges to $T_{u}$. \\ This  
 means that for every $m < \omega$ 
 for every large enough $k \in A^*$ 
  we have  
 $T_{u \cup  \{ k \} }  \cap  {}^{m>} 2 =$ 
 $T_{u} \cap  {}^{m>} 2 .$

\itm  
\begin{enumerate}
\itema $ (\alpha)$ for $A$ an infinite 
  subset of $A^*$  
 we let 
  $$T_{A}  = \bigcup \{ T_{ A \cap n } : n<\omega \} 
\qquad \qquad T[A]  = \bigcup \{ T_{u} : u \subseteq A 
\mbox{ finite}  \}$$ 
\itema  $(\beta)$ for $u,v$  finite 
  subsets of $A^*$ 
 we let  $  {\bf n}(u,v )$ the smallest $m$ such that 
\begin{itemize}
\item whenever  $\eta,\nu$  
 are distinct members of   
 $\lim (T_{u}) \cup  \lim(T_{v})$  
 then the length of $\eta \cap \nu$ 
  is $< m$  
\item $\sup  (w^* \cup   u  \cup  v ) < m $ 
\end{itemize}
Note that ${\bf n}(u,v)$ is well 
  defined,  as both  
 $\lim (T_{u})$ and $\lim (T_{v})$  
 are finite.
\itema  
 $(\gamma)$ for $m< \omega$ let  
 $${\bf n}(m) =^{df} \max \{ {\bf n}(u,v) : u,v \subseteq A^* \cap
(m+1) \},$$  
 so ${\bf n}(m) <\omega$ is well defined  
 being the maximum of a finite 
  set of natural numbers. 
\end{enumerate}

\itm \LABEL{forces} Note that $(w^*, A) \forces \name T \subseteq T[A]$.

 So without loss of generality (as we can  
 replace $A^*$ by any 
  infinite subset):

\itm  if $n\in A^*$ then  
 ${\bf n}(n) < \min  (A^* \setminus (n+1))$

 hence  
  
\itm  if $n \in A^*, u  \subseteq A^*$  
 $\cap (n+1), k \in A^* \setminus (n+1)$  \\
 then $T_{u} \cap {}^{{\bf n}(n)} 2 =$  
 $T_{u \cup \{k\}} \cap$  
  {}$^{{\bf n}(n)} 2 .$

\itm  if $u,v$ are finite subsets  
 of $A^*$ and $\sup  (u \cup v ) < m \in A^*$   
\\  then 
 $T_{u \cup \{m\}} \cap T_{v}  \subseteq$   
 $T_{u} \cap T_{v}$

 Hence  
  
\itm   if $u, v$ are finite  
 subsets of $A^*$, 
 not disjoint for notational  
 simplicity, and  
 $m = \sup( u \cap v )$ then   
$$T_{u} \cap T_{v} \subseteq   T_{u\cap (m+1)} \cap T_{v\cap (m+1)}$$

% $$T_{u} \cap T_{v} \subseteq  \bigcup  
%  \{ T_{u_1} \cap T_{v_1} :\  u_1 \subseteq u \cap (m+1), 
% v_1 \subseteq v \cap (m+1) \}$$
  [ why? we can prove this  by induction  
  on max$(u \cup v )$  
 using (*)$_8 ]$

 Hence, letting $Y =\bigcup \{\lim (T_{u}) :$  
 $u$ a finite subset of $A^*\}$ (=a countable set), we have  
  
\itm  \LABEL{disjoint} if $A, B$ are infinite  
 subsets of  
 $A^*$, with intersection finite  
 non empty\\ then 
 $\lim T[A] \cap\lim  T[B]$ is included in 
$$ \bigcup \{\lim  T_{u} \cap \lim T_{v} :  u \subseteq 
A \cap (\max(A \cap B ) + 1), \  
 v \subseteq 
 B \cap (\max(A \cap B ) + 1) \}$$  
so  $\lim T[A] \cap\lim  T[B]$ is 
is a finite subset of  $Y$.

 [ why ? just use (*)$_9$ and the  
 definition of $T[A] ].$

Now fix an uncountable family $(A_i:i<2^{\aleph_0} )$ of almost disjoint
subsets of $A^*$.  Let $q_i:= (w^*, A_i)$. Then by  \ref{forces}, we
have 
$q_i \forces \name T \subseteq T[A_i]$.   Hence $(q_i, A_i: i < 
2^{\aleph_0} )$ satisfies the
 conditions 1.,2.,3.,4.\ of lemma \ref{mainlemma}. 
\end{proof}

\section{An iteration argument}

\begin{Notation}
For any $f: \omega \to \omega $ we define 
\begin{itemize}
\item $ f^-: \omega \to \omega $  by: $f^-(n)=f(n-1) $ for $n>0$,
$f^-(0)=1$.
\item $ \hat f: \omega \to \omega $  by: $\hat f(n)=2^{f(n)} $.
\end{itemize}
\end{Notation}

\begin{Framework} \label{frame}
We will start with a universe where $2^{\aleph_0} = {\aleph_1}$, 
$2^{\aleph_1} = {\aleph_2}$ and use an iteration of length 
$\kappa = {\aleph_2} $,  
 $\bar{{\Bbb Q}} = \langle {\Bbb P}_{i} , {\Bbb Q}_{i} : i<{\aleph_2}
 \rangle$   satisfying the following: 
\begin{enumerate}
% \item[(A)]  $\kappa$ is a regular cardinal such that for every  
%   $\mu <\kappa$ we have $\mu^{\aleph_0} <\kappa$ 
\item[(A)]  $\bar{{\Bbb Q}} $ 
  is a countable support iteration of proper forcing notions.  
\item[(B)]  $  {\Bbb P}_{\kappa}$,  the union of ${\Bbb P}_{i}$ for
$i<\kappa$ ,    
   satisfies the $\kappa-$cc
\item [(C)]
The set 
$$S:= \{\delta  < \kappa : cf(\delta )>{\aleph_0} , \ \forces_{\P_i} 
\mbox{``$\Q_\delta    $=Mathias forcing, with generic real $\name
\eta_\delta $''}\}$$
is stationary.

\item [(D)]  Each  forcing notion $\name{\Bbb Q}_{i}$  
has the Laver property. 
% is (forced to be) 
%  $(\hat{\name\eta}_\delta, \name\eta^-_\delta ) $-bounding. 
\item[(E)]  Whenever $\alpha < \kappa $ and $\name T$ is a
 ${\Bbb P}_{\alpha}-$name of an Aronszajn 
 tree, then for some $i\in [\alpha ,  \kappa )$,  
 $\Vdash_{{\Bbb P}_{i}}$ ``if $\name T$ is a Souslin  
 tree then $\name {\Bbb Q}_{i}$ is forcing 
 by $\name T$ (or an isomorphic  
 forcing notion)'' 
  
\end{enumerate}
(See the remarks \ref{remlaver} and \ref{remkappa} for weaker assumptions)
\end{Framework}

\begin{Fact}  Let $\bar{{\Bbb Q}}$ satisfy properties (A)--(E)  above,
and let $\P_\kappa $ be the CS limit of this iteration. 
Then:
\begin{enumerate}
\item ${\Bbb P}_{\kappa}$ is proper, making $\kappa$ to $\aleph_2$. 
%  with dense subset of cardinality 
%  $\kappa$ if each ${\Bbb Q}_{i}$ is forced to have  
%  cardinality $<\kappa$ 
%  (or is of cardinality at most $\kappa$ and satisfies the
%  $\kappa-$pic , see [Sh:$f$,Ch  
%  VIII ]. 
\item   If  $\delta \in S$  then  
  the forcing notion  
  ${\Bbb P}_{\kappa} / {\Bbb P}_{\delta+1}$
has the Laver property. 
%  is 
% $(\hat{\name{\eta}}_{\delta} ,  \hat{\name{\eta}}_{\delta}^-)$-bounding. 
\end{enumerate}
\end{Fact}

\begin{proof}  (1) By \cite{Sh:f} 

 (2)By \cite[ch VI, section 3]{Sh:f} 
\end{proof}

\begin{Remark} 
Assume  (say) GCH, then there is a forcing iteration as above. 
Define $\Q_i$
as follows:   If $i$ is even, then let $\Q_i$ be the Mathias forcing, and 
if $i$ is odd, then let $\Q_i$ be either trivial or a Souslin tree. 

Note that in all intermediate universes we will have GCH, all forcing
notions $\P_i$ and $\Q_i$ (for $i<{\aleph_2}  $) will have a dense
subset of size ${\aleph_1}$; this will be sufficient for $\kappa$-cc. 
  All the forcing
notions  $\Q_i$ will have the Laver property (recall that a
Souslin tree does not add reals), and the usual bookkeeping argument
can take care of killing all Souslin trees on $\omega_1$. 
\end{Remark}

\begin{Theorem}\label{mainthm}
 Let $\bar{{\Bbb Q}}$ satisfy conditions (A)--(E) above.  

Then in the universe  ${\bold V}^{{\Bbb P}_{\kappa}}$  
  the following holds:
\begin{enumerate}
\item  Souslin's hypothesis
\item Any nontrivial ccc forcing notion adds a real. (See \ref{known} below)
\item   Whenever $\T$ is a function satisfying the following conditions 
$(\alpha)$--$(\gamma)$:
\begin{itemize}
\item [$(\alpha )$]
	 $\Dom({\Bbb T} )  = \{ f : f$ is a function 
   from $\omega$ to $\omega$, (strictly) increasing $\}$ 
\item [$(\beta )$]
$  {\Bbb T} (f)$ is a subtree of {}$^{\omega>} 2$ 
\item [$(\gamma)$]
 for every $f\in \dom({\Bbb T} )$ and $n<\omega$  
  we have 
     $$1 \le | {\Bbb T}(f) \cap {}^{f(n)} 2|  \le f^-(n )$$ 
\end{itemize}
then $\T$ also satisfies condition $(\delta)$: 
\begin{itemize}
\item[$(\delta)$]
there is a countable subset $Y \subseteq    {}^{\omega} 2$ and an 
 an uncountable subset ${\bf A }\subseteq
   \Dom(\T)$ with: 
\begin{quote}
   whenever  $f \not= g$ are from $A$ then   
   $\lim({\Bbb T}(f )) \cap \lim({\Bbb T}(g ) )$ 
   in a finite subset of $Y.$ 
\end{quote}
\end{itemize}
\end{enumerate}
\end{Theorem}

 \begin{Fact}  \label{known}
 If there is a nontrivial ccc forcing which does not add reals, then
 there is a Souslin tree on $\omega_1$. 
 
 In other words:  If 
  Souslin Hypothesis holds then 
 \begin{quote}
  for  every ccc forcing ${\Bbb R} $   
   which is not trivial, there are  
  $p\in {\Bbb R}$ and an $ {\Bbb R}$-name  
  {$\name \eta$} such that: \\
  $p \Vdash_{{\Bbb R}}$ ``${\name \eta} \in {{}^ {\omega}2}$  
  is new, that is does not 
   belong to ${\bf V}$''
 \end{quote}
 \end{Fact}

\begin{proof}
  Let $\R$ be a nontrivial ccc forcing.   
  So for some $q\in {\Bbb R}$ we have  
\begin{quote}   
  $q \Vdash_{{\Bbb R}} $``$ \name{G}_{{\Bbb R}}$ does  
  not belongs to {\bf V}'' 
\end{quote}
  Hence for some quadruple 
  $(p,\alpha,\beta,{\name \eta})$  
  we have: 

\itm  $p\in {\Bbb R}, \alpha,\beta$ are ordinals,  
 {$\name \eta$} is a ${\Bbb R}$-name 
 and $p\Vdash_{{\Bbb R}} $``$ {\name \eta} \in {}^{\alpha} \beta$  
 is not from {\bf V}'' 
  
 We can choose such quadruple 
  with  
 the ordinal $\alpha$ minimal.   
 Necessarily $\alpha$ is a limit 
  ordinal and for $\gamma<\alpha$,  
 {$\name \eta$} $\restriction \gamma$ 
 is forced by $p$ to belong to {\bf V}. 
  
 For $\gamma<\alpha$, let  
  $$T_{\gamma} = \{ \nu : \mbox{ $\nu$ is a function  
 from   $\gamma$ to $\beta$ so $p$ does not force 
  that $\nu \not= {\name \eta} \restriction \gamma$} \}$$
 and let $T = \cup\{ T_{\gamma} :\gamma<\alpha\}$.  
 Clearly $T$ is a tree with $\alpha$ levels,  
 and  $p$ forces that  
  {$\name \eta$} is a new $\alpha-$branch 
  of it. 
  
 Now $T$ cannot have $\aleph_1$  
 pairwise incomparable  
 elements, as if $\nu_{\zeta}$ for $\zeta<\omega_1$ 
 are like that, we can  
 find $p_{\zeta}$ such that: 
 $p \le p_{\zeta} \in {\Bbb R}$ and 
  $p_{\zeta} \Vdash_{{\Bbb R}} $``$ \nu_{\zeta}$ is an 
  initial segment of {$\name \eta$}'';  
 now if $p_{\zeta},p_{\xi}$ are compatible in  
 ${\Bbb R}$ then $\nu_{\zeta} , \nu_{\xi}$ are 
  comparable  
 in $T$ (being, both, the initial segment  
 of some possible {$\name \eta$}).   
  So $\{ p_{\zeta}:\zeta<\omega_1\}$ are  
 pairwise incompatible 
 contradiction to ${\Bbb R}$ satisfies 
 the ccc'' 
  
 Also in $T$, by its choice, every 
  member has above it elements 
  of every higher level 
 and there is no node above which 
  the tree has no two distinct  
  members of the same level  
 (as then {$\name \eta$} will be  
 forced to belongs to {\bf V}  
 by some 
 condition above $p)$. 
  
 Also as ${\Bbb R}$ satisfies the ccc, 
  every level is countable,  
 and by the minimality of  
 $\alpha$ (as we are allowed to change $\beta)$ 
 clearly $\alpha$ is a regular  
 cardinal.  
 Now $\alpha> \omega_1$ 
 is impossible by 
  ``${\Bbb R}$ satisfies the $\kappa-$cc''. 
 As there are no Souslin 
  tree also $\alpha=\omega_1$ is  
 impossible. So clearly $p$ 
  forces that ${\Bbb R}$ add reals  
 so there is {$\name \eta$} as  
 required. 
  
\end{proof}

\begin{Observation}
Theorem \ref{mainthm}  suffices to prove \ref{jlx} and its
  strengthening $(**)$ mentioned in the introduction. 
  That is, conditions (2)\&(3) of \ref{mainthm}
 imply:
\begin{quote}
Any nontrivial forcing with the Laver property has an uncountable 
antichain. 
\end{quote}
\end{Observation}

\begin{proof}
  Let $\R$ be a forcing notion with the Laver
property which adds a real, say $p \forces \name \eta \in {}^\omega 2
, \name \eta\notin V$. 

Consider any increasing function $f$. The function $n\mapsto 
\eta\on f(n+1)$ has only $2^{f(n+1)}$ many possible values, i.e., is
bounded.  So, by the Laver property there is a tree $T_f \subseteq
\tlo$ and a condition $q_f$ stronger than $q$  with
$$ \forall n\ | T_f \cap 2^{f(n+1)}| \le f(n), 
\qquad \qquad q_f \forces \name \eta \in \lim T_f $$
We have thus defined a family $\T=(T_f: f
  \mbox{ increasing })$.
  By theorem
\ref{mainthm}, there is a family 
$(f_i:i\in \omega_1)$ such that 
\begin{quote}
$ \forall i \not= j: $  $\lim T_{f_i} \cap \lim T_{f_j}$ is finite
\end{quote}
Clearly, for $i\not=j$ the conditions $q_{f_i}$  and  $q_{f_j}$  must
be incompatible, since any condition $ r $ stronger than both would 
force 
$$ r \forces_\R \name \eta \in \lim T_{f_i}  \cap  \lim T_{f_i} $$
which implies $r\forces \name \eta\in V$, a contradiction. 

(Remark: While $\lim T_{f_i}$ and $\lim T_{f_j}$ can of course contain 
branches  in $V^\R$ which did not exist in $V$, the fact that their
intersection is a certain  finite set is absolute between $V$ and $ V^\R$.)
\end{proof}

\subsection*{Proof of theorem \ref{mainthm}, part 3}

Assume that
 {$\name {\Bbb T}$} is a
$\P_\kappa$-name      such that 
 $p^* \in \P_\kappa  $  forces ``$\name \T$ satisfies
 $(\alpha),(\beta),(\gamma) $''.

Without loss of generality (replacing the ground model by an
intermediate model $V[G_\alpha]$, $G_\alpha \subseteq \P_\alpha$,
$p^*\in G_\alpha$, if necessary) we can assume that $p^* $ is really
the empty condition.  

Let $S \subseteq \kappa  $ be unbounded, $\delta \in S \Rightarrow 
\Q_\delta$ is Mathias forcing, with generic real $\name\eta_\delta$. 

 For every $\delta\in S$ let $\name T^0_{\delta} =  
  \name {\Bbb T} (\name \eta^0_{\delta} )$;  clearly : 

\itm $   \Vdash_{{\Bbb P}_{\kappa}} $``$ {\name T}^0_{\delta}$  
  is a subtree of ${}^{\omega>} 2$ such that  
   $(\forall n ) |{\name T}^0_{\delta} \cap
        {}^{{\name \eta}_{\delta} (n)}2 | 
 \le  {\name \eta}^-_{\delta} (n)$.''

There is only a bounded number of
 possibilities for 
   ${\name T}^0_{\delta} \cap
  {}^{{\name \eta}_{\delta} (n)}2$, so since  
${\Bbb P}_{\kappa} / {\Bbb P}_{\delta+1}$ has the Laver
  property, we can find a pair  
 $(p_{\delta} , {\name T}_{\delta})$ 
    satisfying 
   \itm  
\begin{enumerate} 
\item   $p_{\delta}\in \P_\kappa $  
\item   
 $ {\name T}_{\delta}$ is a  
  ${\Bbb P}_{\delta+1}-$name   
\item 
 $  \Vdash_{{\Bbb P}_{\delta+1}}  $``$   {\name T}_{\delta}$ 
   is a subtree of {}$^{\omega>}2$   
   and $n<\omega \rightarrow$ 
   $\vert {\name T}_{\delta} \cap {}^{{\name \eta}_{\delta} ( n )} 2 \vert$ 
  $\le \hat {\name \eta}^-_{\delta} (n) $'' 
\item    $  \Vdash_{{\Bbb P}}   $``$   {\name T}^0_{\delta}$  
   $\subseteq {\name T}_{\delta}  $'' 
\end{enumerate}

  So there is  stationary subset $S_1 \subseteq S$ 
 and  a condition $q_1\in \P_\kappa $ such that 
  $\delta \in S_1 \rightarrow p_{\delta} \restriction \delta = q_1 .$  
(Again we may assume that $q_1$ is the trivial condition.)

 Possibly increasing $p_{\delta} (\delta )$  
 we can find a $\P_\delta $-name such that  
$p_\delta \on \delta $ forces: 
\itm 
Above $p_\delta(\delta)$, 
the  $\Q_\delta $-name ${\name T}_{\delta}$ can  
  be read continuously  
  from {$\name \eta$}$_{\delta}$ as in  \ref{continuous}, through the 
function $\name t_\delta $).

For $\delta\in S_1 , p_{\delta} (\delta)$,   $\name t_\delta$
 {$\name T$}$_{\delta}$ are members of  
  ${\mathcal H}(\aleph_1 )^{{\bold V} [ G_{\delta} ]}$.  So we can
 find $q_\delta \ge p_\delta\on \delta $ forcing
  $p_\delta(\delta)$,
 $\name t_\delta$ and  $\name T_\delta$ to be equal to hereditarily
 countable $\P_\delta$-names 
  $p'_\delta(\delta)$,
 $\name t'_\delta$ and  $\name T'_\delta$. 
[Here, ``hereditarily countable'' is taken in the sense of 
 \cite[III 4.1A]{Sh:f}.]

Since $cf(\delta)>{\aleph_0}$ for $\delta\in S_1$ we can find a
  stationary subset $S_2 \subseteq S_1$ on which   $p'_\delta(\delta)$,
 $\name t'_\delta$ and  $\name T'_\delta$ are all constant, say with
  values    $p^*$,
 $\name t^*$ and  $\name T^*$. Again we change our base universe to some
  intermediate universe so that
  $\name T^*$ is now a $\M$-name, and $t^*$ is an actual function, and
  $p^*= (w^*,A^*)\in \M$. 

We now use our main lemma \ref{mainlemma} to find an almost disjoint
family $(A_i: i \in \omega_1)$ and  $(T_i: i \in \omega_1)$ such that $A_i
\subseteq A ^*$, $(w^*,A_i) \forces_\M \name T^* \subseteq T[A_i]$.

Note because of \ref{continuous} this relation can already be computed
from $t^*$, so we also have: 

\itm \LABEL{sixteen} $ \forall \delta \in S_2 $ 
  $ \forces_{\P_\delta}$``$  q_i \forces _{\Q_\delta} \name  T_\delta
  \subseteq T[A_i] $''. 

Now consider the  model $V[G_\kappa]$.  A density argument shows that 
\itm  $ \forall i: \{ \delta \in S_1: q_i\in G_{\Q_\delta}\} \not=
\emptyset$

So for all $i$ there is $\delta = \delta(i)$ with 
$\T(\eta_{\delta(i)}) \subseteq T[A_i]$.   Letting 
$${\bf  A}:= \{ \delta_i: i < \omega_1 \}$$
we have found an uncountable family as required.

\section{Refinements} 

Theorem \ref{mainthm} answers the original question, but 
 essentially the same proof gives a somewhat
stronger theorem.
 The following remarks point a few places where
assumptions can be weakened or conclusions strengthened.  We leave
the details to the reader.

\begin{Remark} \label{remkappa}
$2^{2^{\aleph_0}} = {\aleph_2}$ in the ground model is
 not necessary.  The length of our iteration can be any
regular cardinal $ \kappa $ satisfying $ \mu^ {\aleph_0} < \kappa $
for all $\mu < \kappa $.    In the final model we will  have
$2^{\aleph_0} = {\aleph_2} = \kappa $. 

\end{Remark}

\begin{Remark} \label{remlaver}
  It is not necessary that all forcing notions have the
Laver property.   All we need is that $\P_\kappa / P_{\delta+1}$ is 
$(\hat{\name \eta}_\delta , \hat{\name \eta}^-_\delta / \name
\eta_\delta ^-)$-bounding, which can be ensured by a slightly stronger
condition on all the $\Q_i$, $\delta < i < \kappa $. 
\end{Remark}

\begin{Remark} 
We showed that in our model every forcing notion with the
Laver property which adds reals will have an uncountable antichain. We
can strengthen this conclusion by remarking that such forcing notions
will actually have an antichain of size $\kappa = {\aleph_2}$. 
\end{Remark}
\begin{proof}
Recall the construction of the almost disjoint family after
condition \ref{disjoint}, which was used in \ref{sixteen}. 
Instead of using an almost disjoint family of size ${\aleph_1}$ in the
intermediate model we can use a $\P_\kappa$-name of an
almost disjoint family of size continuum: Identify the set $A^*$ there
with $\tlo$, then every $\P_\alpha$-name $\name\rho$ of an element of 
${}^\omega 2$ will induce a set $A_{\name \rho} \subseteq A^ *$.  
Clearly, 
$$ \forces_{\P_\kappa } \name \rho_1 \not= \name\rho_2 \Rightarrow 
A_{\name \rho_1} \cap A_{\name \rho_2 }  \mbox{ finite}. $$
As before, a density argument ensures that there will be 
$\kappa$ many different functions $\name \rho$ such that  
 $(w^*, A_{\name \rho})$ appears in one of the generic Mathias filters
for some $\Q_\delta=\M$, so 
we can strengthen the conclusion in \ref{mainthm}, 3$(\delta)$ to get
a $\kappa$-size set $\bf A$  rather than just an uncountable one. 
\end{proof}

\begin{Remark}  
We  do not need that all forcing notions $\Q_i$ have size at most
${\aleph_1}$, there are also weaker conditions (e.g. $\kappa$-pic, see
\cite[Ch VIII]{Sh:f}) that will ensure $\kappa$-cc of $\P_\kappa$. 

For example, instead of forcing only with Souslin trees in the odd
 stages we can   use 
 the forcing from 
 \cite[Ch V, Section  6]{Sh:f}, 
 it specializes the tree (so we  
 can specialize all  
 Aronszajn trees).   Here we can prove the $\kappa$-cc 
 using the $\kappa$-pic condition. , 
\end{Remark}

\begin{Remark}
Finally, in $V^{\P_\kappa }$ we can strengthen the conclusion 
\begin{quote} Every ccc nontrivial forcing fails the Laver property
\end{quote}
as follows: 
\begin{quote}
For every ccc forcing notion  ${\Bbb R}$, whenever 
 $\name\eta$ is an {${\Bbb R}$}-name  
 of a new member of   ${}^ \omega  2$ 
 and $h$ is a strictly increasing  
 function from  
 $\omega$ to $\omega, \\ \bf{{\rm then}}$ we can  
 find an increasing sequence 
 $\langle n_{i}:i<\omega\rangle$ of natural 
  numbers such that $h(n_{i}) < n_{i+1}$  
 and  for no $p,T$ do we have: 
\begin{quote}
 $p\in {\Bbb R}, T$ a subtree  
 of $\tlo $,  
 $p \Vdash_{{\Bbb R}} $``$ \name \eta \in 
  \lim(T)$ and for every $i<\omega$ 
 we have $\vert   T \cap {}^{n_{2i+1}} 2 \vert$    
 $\le h(n_{2i} )$ 
\end{quote}
\end{quote}
\end{Remark}
The proof is similar to the proof above.

  \bibliographystyle{lit-unsrt}
\bibliography{listb,listx,lista}

\end{document}